\documentclass{amsart}

\makeindex

\usepackage[all]{xy}
\usepackage[english]{babel}

\usepackage{amsmath}

\usepackage{amsthm}

\usepackage{amscd}

\usepackage{amsfonts}

\usepackage{amssymb}

\usepackage{mathrsfs}

\newtheorem{Lem}{Lemma}[section]

    \newtheorem{Prop}[Lem]{Proposition}

    \newtheorem{Thm}[Lem]{Theorem}

    \newtheorem{Cor}[Lem]{Corollary}

    \newtheorem{lemma}{Lemma}

\theoremstyle{definition}

    \newtheorem{Rem}[Lem]{Remark}

\newcommand{\Ps}{\mathbb{P}^}
\newcommand{\Z}{\mathbb{Z}}

\newcommand{\ra}{\rightarrow}
\newcommand{\lra}{\longrightarrow}
\newcommand{\R}{\mathbb R}

\newcommand{\C}{\mathbb C}

\begin{document}
\title{A compactification of $\mathcal M_3$ via $K3$ surfaces}
\author{Michela Artebani}
\address{Dipartimento di Matematica, Universit\`a di Milano, via C. Saldini 50, Milano, Italia}
\email{michela.artebani@unimi.it} \subjclass[2000]{14J10, 14J28, 14H10} \keywords{genus three curves, K3 surfaces}
\thanks{This work is partially supported by: PRIN 2003: Spazi di moduli e teoria di Lie; GNSAGA}
\begin{abstract}
S. Kond\=o defined a birational period map from the moduli space of genus three curves to a moduli space of degree four polarized K3 surfaces. In this paper we extend the period map to a surjective morphism on a suitable compactification of $\mathcal M_3$ and describe its geometry.
\end{abstract}
\maketitle
\section*{Introduction}
Let $V=\mid\!\mathcal O_{\Ps 2}(4)\!\mid\cong\Ps{14}$ be the space of plane quartics and $V_0$ be the open subvariety of smooth curves.
The degree four cyclic cover of the plane branched along a curve $C\in V_0$ is a $K3$ surface equipped with an order four non-symplectic automorphism group.
This construction defines a holomorphic period map:
$$\mathcal P_0:\mathcal V_0\longrightarrow \mathcal M,$$
where $\mathcal V_0$ is the geometric quotient of $V_0$ by the
action of $PGL(3)$ and $\mathcal M$ is a moduli space of polarized
$K3$ surfaces.


In \cite{K} S. Kond\=o shows that $\mathcal P_0$ gives an
isomorphism between $\mathcal V_0$ and the complement of two
irreducible divisors $\mathcal D_n, \mathcal D_h$ in $\mathcal M$.
Moreover, he proves that the generic points in $\mathcal D_n$ and
$\mathcal D_h$ correspond to plane quartics with a node and to
smooth hyperelliptic genus three curves respectively.

The moduli space $\mathcal M$ is an arithmetic quotient of a six
dimensional complex ball, hence a natural compactification is given
by the Baily-Borel compactification $\mathcal M^*$ (see \cite{BB}).
On the other hand, geometric invariant theory provides a compact
projective variety $\mathcal V$ containing $\mathcal V_0$ as a dense
subset, given by the categorical quotient of the semistable locus in
$V$ for the natural action of $PGL(3)$. In this paper we prove that
the map $\mathcal P_0$ can be extended to a holomorphic surjective
map
$$\mathcal P:\widetilde{\mathcal V}\lra \mathcal M^{*}$$
on the blowing-up  $\widetilde{\mathcal V}$ of $\mathcal V$ in the point $v_0$
corresponding to the orbit of double conics. In particular, we
show that the exceptional divisor in $\widetilde{\mathcal V}$ is a GIT moduli space of hyperelliptic genus three curves which is mapped isomorphically onto the Baily-Borel compactification of $\mathcal D_h$.

The paper is divided in three main sections related to plane
quartics, hyperelliptic genus three curves and finally genus three
curves in general.\\

The first section starts reviewing the geometric construction by
Kond\=o, which defines the period map  $\mathcal P_0$ on the moduli
space of smooth plane quartics.
This construction can be easily extended to stable
quartics (i.e. having at most ordinary nodes and cusps) since in
this case the 4:1 cyclic cover of the plane branched along the curve
has at most rational double points, hence its minimal resolution is
still a K3 surface. In fact the period map can be extended to a
holomorphic map
$$\mathcal P_1:\mathcal V_s\lra \mathcal M$$
defined on the geometric quotient $\mathcal V_s$ of the locus of
stable quartics for the action of $PGL(3)$. The geometry of $K3$
surfaces associated to stable quartics is also described in detail.

The analogous construction for strictly semistable quartics gives
surfaces with elliptic singularities for quartics with tacnodes and
surfaces with significant limit singularities (see \cite{S3}) in the
case of double conics. Nevertheless, by taking the Baily-Borel
compactification $\mathcal M^*$ of the ball quotient $\mathcal M$,
we can extend the period map to a morphism
$$\mathcal P_2:\mathcal V_{ss}\backslash\{v_0\}\lra \mathcal M^*,$$
where the point $v_0$ represents the orbit of double conics and
strictly semistable quartics give a smooth rational curve (with
$v_0$ in its closure) mapped to the boundary.

In section two, we give a correspondence between hyperelliptic genus three curves and certain hyperelliptic polarized $K3$ surfaces parametrized by $\mathcal D_h$.
In fact, this relation has been studied in detail by Kond\=o in \cite{Kh}.
He defines a period map
$$\mathcal P^h:\mathcal V^h\lra \mathcal D^*_h$$
from the GIT moduli space $\mathcal V^h$ of sets of eight points in
$\Ps 1$ to the Baily-Borel compactification of $\mathcal D_h$ and
proves that it is an isomorphism. We recall his results giving a
different approach. In particular, we show that $\mathcal D_h$ is a
moduli space for degree four cyclic covers of a cone in $\Ps 5$
branched along a quadratic section.

The last section contains the main theorem. We first construct a
blowing up $\widetilde{\mathcal V}$ of $\mathcal V$ in $v_0$ such
that the exceptional divisor is isomorphic to the moduli space
$\mathcal V^h$. We then prove that the period maps $\mathcal P_2$
and $\mathcal P^h$ define a global period map
$$\mathcal P:\widetilde{\mathcal V}\lra \mathcal M^*.$$
Let $\widetilde{\mathcal V}_s$ be the locus corresponding to stable
genus three curves in $\widetilde{\mathcal V}$. The morphism
$\mathcal P$ induces an isomorphism $\widetilde{\mathcal V}_s\cong
\mathcal M$ and maps the smooth rational curve $\widetilde{\mathcal
V}\backslash\widetilde{\mathcal V}_{s}$ to the boundary of $\mathcal
M^*$ (which consists of one point).
\\

Similar descriptions for sextic double planes have been given by E.
Horikawa in \cite{H1} and by J. Shah in \cite{S2}. These two papers
and \cite{S} by H.J.M. Sterk are all important references for
this work.\\

The Appendix contains a brief review on GIT moduli spaces of genus three curves.\\

\emph{Acknowledgements.} This paper is part of my PhD thesis. I'm
grateful to my advisor, Prof. B. van Geemen, for helpful comments
and careful reading. I also would like to thank Dr. \!A. Laface and
Prof. E. Looijenga for several interesting discussions and
suggestions. I acknowledge the Mathematics Department of the
University of Milano for supporting me.
\section{Plane quartics}
\subsection{Smooth quartics}
We start recalling the geometric construction introduced by Kond\=o
in \cite{K}. This can be resumed by the following commutative diagram:
$$\xymatrix{
X_C \ar[rd]_{\pi_2} \ar@{-}[rr]^ {\pi}& \ar[r]&  \Ps 2\\
  &S_C \ar[ur]_{\pi_1}&  },$$
where $C$ is a smooth quartic curve, $\pi_1$ is the
double cover of the plane branched along $C$ and $\pi$
is the 4:1 cyclic cover of the plane branched along $C$. Note that
$S_C$ is a Del Pezzo surface of degree two and $\pi_1$ is the
morphism associated to the anti-canonical linear system. Its double
cover $X_C$ is a $K3$ surface. In coordinates, if the quartic $C$ is
defined by the equation $f_4(x,y,z)=0$:
$$X_C=\{(x,y,z,t)\in \Ps 3:\ t^4=f(x,y,z)\}.$$
The $K3$ surface $X_C$ has a natural degree four polarization
induced by the embedding in $\Ps 3$ and given by $\pi^*(l)$, where
$l$ is the class of a line in the plane.

Let $G$ be the covering transformation group of $\pi$. A generator
$\sigma$ for $G$ can be chosen such that the space of holomorphic
two-forms of $X_C$ lies in the $i$-eigenspace $W$ of the isometry
$\rho=\sigma^*$ on $H^2(X_C,\C)$. In particular the Picard lattice of $X_C$ contains the pull-back $L_+$ of the Picard lattice of $S_C$ as a sublattice.
Hence the period point of $X_C$ lies in the six dimensional complex ball
$$B=\{x\in \mathbb P(W): (x,\bar x)>0\}\subset \mathbb P(L_-\otimes \C),$$
where $L_-$ is the orthogonal complement of $L_+$ in $H^2(X_C,\Z)$.

Following the definition given in \cite{DVK}, an
$(L_+,\rho)$-\emph{polarized} $K3$ surface $X$ is an $L_+$-polarized
$K3$ surface with period point in $B$ (up to the choice of an
isometry $H^2(X,\Z)\cong H^2(X_C,\Z)$).  As proved in \cite{K}, the
moduli space $\mathcal M$ of these polarized $K3$ surfaces is the
quotient of the six dimensional complex ball $B$ by the action of an
arithmetic group $\Gamma$ of automorphisms:
$$\mathcal M=B/\Gamma,\ \Gamma=\{\gamma\in O(L_-): \gamma\circ \rho=\rho\circ\gamma\}.$$
Projectively equivalent plane quartics give isomorphic polarized $K3$ surfaces, thus the above construction defines a holomorphic period map:
$$\mathcal P_0:\mathcal V_0\lra \mathcal M,$$
where $\mathcal V_0\cong\mathcal M_3\backslash\mathcal M_3^h,$ is
the moduli space of  smooth quartics (see the Appendix).

Note that all $(L_+,\rho)$-polarized $K3$ surfaces have a degree four
polarization given by the $\rho$-invariant lattice $L^{\rho}=\langle h\rangle$ in $L_+$. The polarization is \emph{ample} if there are no
$(-2)$-curves orthogonal to $h$, i.e. a $K3$ surface is not ample polarized iff
its class belongs to the \emph{discriminant locus} $\mathcal
D$ in $\mathcal M$:
$$\mathcal D=\left(\bigcup_{r\in \Delta}H_{r}\right)/\Gamma,$$
where $\Delta=\{r\in L_-:\ r^2=-2\}$ is the set of \emph{roots} of
$L_-$ and $H_{r}=B\cap r^{\perp}.$ In fact, ample $(L_+,\rho)$-polarized $K3$ surfaces correspond to smooth plane quartics.
\begin{Thm}[Theorem 2.5, \cite{K}]
The period map gives an isomorphism:
$$\mathcal P_0:\mathcal V_0\cong \mathcal M_3\backslash\mathcal M_3^h\lra \mathcal M\backslash \mathcal D.$$
 \end{Thm}
If $\Lambda_r=\langle r,\rho(r)\rangle$ and $\Lambda^{\perp}_r$ is its orthogonal lattice in $L_-$, then the roots can be divided in two classes (Lemma 3.3, \cite{K}):
$$\Delta_n=\{r\in \Delta:\ \Lambda_r^{\perp}\cong U^{\oplus 2}\oplus A_1^{\oplus 8}\},\
\Delta_h=\{r\in \Delta:\ \Lambda_r^{\perp}\cong U(2)^{\oplus 2}\oplus D_8 \}.$$
This leads to a natural decomposition of $\mathcal D$ in the union of two divisors, called \emph{mirrors}:
$$\mathcal D_n=\left(\bigcup_{r\in \Delta_n} H_r\right)/\Gamma,\ \mathcal D_h=\left(\bigcup_{r\in \Delta_h} H_r\right)/\Gamma.$$
In the following sections we will see that these two divisors have a clear interpretation in terms of genus three curves.
\subsection{Stable quartics}
A plane quartic is \emph{stable} for the action of $PGL_3$ if and
only if it has at most ordinary nodes and cusps (see the Appendix).
We show that the period map has a good behaviour on the stable
locus.
\begin{Lem}\label{res}
The minimal resolution of the degree four cyclic cover of $\Ps 2$ branched along a stable quartic is a $(L_+,\rho)$-polarized $K3$ surface.
\end{Lem}
\proof Let $\pi:Y\lra \Ps 2$ be the degree four cyclic cover branched along a stable plane quartic $C$. Local computations show that
$\pi^{-1}(p)$ is a rational double point of type $A_3$ if $p$ is a node and of type $E_6$ if $p$ is a cusp of $C$ (see Ch.III, \cite{BPV}).
The minimal resolution $r:\tilde Y\lra Y$ of $Y$ is a $K3$ surface
since $r^*(K_Y)=K_{\tilde Y}=0$ and $H^1(\mathcal O_{\tilde
Y})=H^1(\mathcal O_Y)=0$. Since these $K3$ surfaces are
degenerations of ample $(L_+,\rho)$-polarized $K3$ surfaces their
period points belong to $B$.
\qed\\

By the previous Lemma a polarized $K3$ surface $X_C$ can be
associated to any stable plane quartic $C$. This defines a natural
extension $\mathcal P_1$ of $\mathcal P_0$ to the locus $\mathcal
V_s$ representing stable quartics in $\mathcal V$:
$$\mathcal P_1:\mathcal V_s\lra \mathcal M.$$
\begin{Prop}\label{s}
The period map $\mathcal P_1$ is a holomorphic extension of
$\mathcal P_0$ and maps $\mathcal V_s\backslash \mathcal V_0$ to
$\mathcal D_n.$ Moreover, the map is generically surjective onto
$\mathcal D_n$.
\end{Prop}
\proof By a result of E. Brieskorn on simultaneous resolution of
singularities (see \cite{Br}) the period map $\mathcal P_0\circ q$
(see the Appendix) can be extended holomorphically to the open
subset of stable curves in $V$. By taking the quotient for the
action of $PGL(3)$ we get the first assertion. The generic point in
$\mathcal V_s\backslash \mathcal V_0$ corresponds to a plane quartic
with one node. This case is analyzed in detail in \cite{K}, in
particular it is proved that the polarization on $X_C$ is not an
ample since, for general $C\in \mathcal V_s\backslash \mathcal V_0$:
$$Pic(X_C)\cap L_-=\langle r,\rho(r)\rangle,\  r\in \Delta_n.$$
This implies the first assertion, the second one is proved in \cite{K}, \S4.
\qed\\

We can define the \emph{singular type} of a stable quartic to be the pair $(n,c)$ where $n$ is the number of its nodes and $c$ of its cusps. The singular type induces a natural stratification of $\mathcal D_n$, since we can associate to each pair $(n,c)$ ($(n,c=\not=(0,0)$) the variety $\mathcal D_{n,c}$ corresponding to $K3$ surfaces
$X_C$ with $C$ having at least $n$ nodes and $c$ cusps.
We now describe the geometry of the generic $K3$ surface in each stratum.

For any stable quartic $C$ of singular type $(n,c)$, the $K3$
surface $X_C$ carries a natural involution $\tau_{n,c}$ in the
covering group of the 4:1 map $X_C\lra \Ps 2$. Let $\tau^*_{n,c}$ be
the induced isometry on $H^2(X_C,\Z)$.
\begin{Prop}\label{stab}
The Picard lattice of the generic K3 surface $X_C$ in $\mathcal D_{n,c}$ equals the invariant lattice of $\tau^*_{n,c}$. Moreover
$$Pic(X_C)\cap L_-\cong A^{\oplus 2n}_1\oplus D^{\oplus c}_4,$$
in particular $X_C$ has Picard number $\rho(n,c)=8+2n+4c.$
\end{Prop}
\proof The case $n=c=0$ has been described in the previuos section and the case of stable quartics with $c=0$ and $1\leq n\leq 3$ is an easy generalizaton of the
results for $n=1$ given in the proof of Proposition \ref{s}.

We study in more detail the case of a plane quartic $C$ with an ordinary cusp at a point $p$.
The minimal resolution of the 4:1 cyclic cover of $\Ps2$ branched along $C$ can be obtained in the following steps (in order to simplify the notation, at each step the name of a curve is the same as that of its proper transform) :\\
1) Let $b_1:X_1\lra \Ps 2$ be the blow up of $\Ps2$ in $p$ and $E$ be the exceptional divisor. Notice that $C$ is smooth and intersects $E$ in one point with multiplicity $2$.\\
2) Take the double cover $\phi_1:S_C\lra X_1$ branched along $C$. Let $\phi^{-1}(E)=E',E''$.\\
3) Let $b_2:X_2\lra S_C$ be the blow up in the point $E'\cap C=E''\cap C$. Let $F$ be the exceptional divisor.\\
4) Let $b_3:X_3\lra X_2$ be the blow up in the points $C\cap(F\cup E'\cup E'')$. Let $G,H,L$ be the exceptional divisors.\\
5) Let $\phi_2: X_C\lra X_3$ be the double cover branched
along the divisor $C\cup F\cup E'\cup E''$.

The surface $S_C$ is a nodal Del Pezzo surface of degree two (see \cite{DO}), hence the pull-back of its Picard lattice in $X_C$ is isomorphic to $L_+$. In particular $X_C$ is a $(L_+,\rho)$-polarized $K3$ surface with
$$Pic(X_C)\cap L^{\perp}_+=\langle F,G,H,L\rangle\cong D_4.$$
Moreover, $X_C$ carries a natural involution $\tau_{0,1}$ induced by the cover $\phi_2$ and
the invariant lattice $L^{0,1}_+$ of the induced isometry $\tau^*_{0,1}$ is a primitive sublattice of $Pic(X_C)$ of rank $r(L_+)+4=12$.
In fact, for the generic $X_C$ in $\mathcal D_{0,1}$ the two lattices coincide by dimension reasons.

The general case (i.e. $n,c>1$) is an easy generalization of the
previous ones.

\qed
\begin{Cor}\label{irr}
The following table lists the isomorphism classes of the Picard
lattices $Pic(n,c)$ of the generic K3 surfaces in $\mathcal
D_{n,c}$, with $0\leq n,c\leq 3$ (hence the general plane quartic
with singular type $(n,c)$ is irreducible):
$$\begin{array}{lllll}
(n,c)&  Pic(n,c)\\
\ \\
(0,0)& \langle 2\rangle\oplus A_1^{\oplus 7}\\
(1,0)& U\oplus A_1^{\oplus 8}\\
(2,0)& U\oplus A_1^{\oplus 6}\oplus D_4\\
(3,0)& U\oplus A_1^{\oplus 6}\oplus D_6\\
(1,1)& U\oplus A_1^{\oplus 2}\oplus D_4\oplus D_6\\
(1,2)& U\oplus A_1^{\oplus 2}\oplus D_6\oplus E_8\\
(2,1)& U\oplus A_1^{\oplus 2}\oplus D_4\oplus D_8\\
(0,1)& U\oplus A_1^{\oplus 4}\oplus D_6\\
(0,2)& U\oplus A_1^{\oplus 2}\oplus D_4\oplus E_8\\
(0,3)& U\oplus A_1^{\oplus 2}\oplus E_8^{\oplus 2}\\
\end{array}$$
\end{Cor}
\proof We adopt the same notation of the proof of Proposition
\ref{stab}. Since $X_C$ is general and $0\leq n,c\leq 3$, we can
assume $C$ to be an irreducible plane quartic. In this case, the
curve $C$ in $X_C$ is irreducible and each double point decreases
its genus by one. Moreover, the fixed locus of the involution
$\tau_{n,c}$ contains a smooth rational curve for every node in $C$
and $3$ smooth rational curves for every cusp. Hence the fixed locus
of $\tau_{n,c}$ is given by an irreducible curve of genus $3-n-c$
and $n+3c$ smooth rational curves. Note that $L^{n,c}_+$ is a
$2$-elementary lattice, let $\ell_{n,c}$ be the minimal number of
generators of its discriminant group and $\delta_{n,c}$ be the
invariant defined in \cite{N1}. By \cite{N2}, Theorem 4.2.2, we get:
$$\rho_{n,c}+\ell_{n,c}=22-2(3-n-c),\ \ \rho_{n,c}-\ell_{n,c}= 2(n+3c).$$
Hence $\ell_{n,c}=8-2c$. Moreover, if $(n,c)\not=(1,1), (1,2)$, we
have $\delta_{n,c}=1$ by \cite{N2}, Theorem 4.3.2. In these cases,
let $x$ be the class of a $-2$-curve in $L_+$ corresponding to a
bitangent of $C$ not passing through the singular points of $Q$.
Then it is easy to see that $x/2$ belongs to the dual lattice of
$L^{n,c}_+$ and $(x/2)^2=-1/2\not\in \Z$. It follows that
$\delta_{n,c}=1$ also in these two cases. Theorem 4.3.2 in \cite{N2}
also gives that the isomorphism class of the lattice $L^{n,c}_+$  is
determined by the invariants $(\rho_{n,c},\ell_{n,c},\delta_{n,c})$.
Hence it is enough to compute this triple of invariants for the
lattices in the second column (see \cite{N2}, Proposition
3.2.2).\qed
\begin{Rem}
A $K3$ surface associated to a singular stable quartic $C$ carries a
natural elliptic fibration, given by the inverse image of the pencil
of lines through a singular point of $C$. If $C$ is the generic
quartic with one node, then this pencil has $8$ singular fibers of
type $III$ in the sense of Kodaira, corresponding to the two
branches of the node and to the six lines through the node and
tangent to $Q$ in smooth points. If $C$ is the generic quartic with
one ordinary cusp the fibration has $6$ fibers of type $III$ and one
of type $I_0^*$. They correspond respectively to the $6$ tangent
lines passing through the cusp and to the tangent line in the cusp.
\end{Rem}
\subsection{Semistable quartics}
Let $C$ be a plane quartic with an ordinary tacnode at $p$ and $Y$ be the 4:1 cyclic cover of the plane branched along $C$. The local equation of $Y$ over the point $p$ is
analytically isomorphic to:
$$t^4+y^2+x^4=0,$$
so $Y$ has an elliptic singularity of type $\tilde E_7$. This
suggests that there will be no $K3$ surfaces corresponding to these
curves. However, we prove that the period map can be still extended
holomorphically to strictly semistable admissible quartics if we
consider the Baily-Borel compactification of the period domain
$\mathcal M$. In particular, we show that these quartics are mapped
to the boundary and we describe the corresponding degenerations of
$K3$ surfaces.

Let $\mathcal M^*$ be the Baily-Borel compactification of $\mathcal
M$ (see \cite{BB}). It easy to prove that the $i$-eigenspace $W$ of
$\rho$ on $L_-\otimes \C$ is given by
$$W=\{x-i\rho(x): x\in L_-\otimes_{\Z}\R\}\cong L_-\otimes_{\Z}\R.$$
We say that $w\in W$ is \emph{defined over} $\Z$ if $w=x-i\rho(x)$
with $x\in L_-\subset L_-\otimes_{\Z}\R$.
\begin{Lem}\label{bou}
The boundary of  $\mathcal M^*$ is a disjoint union of points
corresponding to isotropic lines in $W$ defined over $\Z$.
\end{Lem}
\proof The period domain of $L_+$-polarized $K3$
surfaces is given by:
$$\mathcal D_+=\{z\in \mathbb P(L_-\otimes \C):\ (z,z)=0,\ (z,\bar
z)>0\}.$$ The rational boundary components of $\mathcal D_+$ are
given by points and curves corresponding respectively to isotropic
lines and planes with generators in $L_-$ (see \cite{BB}). Hence the
rational boundary components of $B$ are given by the intersection of
these components with $W$. Notice that an isotropic line $\mathbb
Cv$, $v^2=0$, $v\in L_-$ can not be contained in $W$, since $W\cap
L_-=\{0\}$. Let $\Lambda$ be an isotropic plane generated by $x,y\in
L_-$. If $\Lambda$ contains an element in $W$, then it is of the
form:
$$ax+by=a_1x+b_1y-i(a_1\rho(x)+b_1\rho(y)),$$
where $a_1,b_1\in \R$. By intersecting both sides with $y$,
this gives that $(\rho(x),y)=0$. Hence:
$$\Lambda=\langle x,y\rangle=\langle x,\rho(x)\rangle.$$
Thus the intersection $\Lambda\cap W$ is a line generated by an
isotropic vector defined over $\Z$:
$$ax+by=(a_1+ib_1)(x-i\rho(x)).$$
Conversely, an isotropic vector in $W$ defined over $\Z$ is of the
form $x-i\rho(x)$, $x\in L_-$ and the plane generated by $x,
\rho(x)$ is an isotropic plane with generators in $L_-$.
\qed\\

Let $q:V_{ss}\lra \mathcal V$ be the quotient morphism to the
categorical quotient of $V_{ss}$ by the action of $PGL_3$ and
$V_{ss}'\subset V$ be the set of admissible semistable plane
quartics (see the Appendix). The following result and its proof are
similar to those of Theorem 4.1, \cite{H1} and Theorem 5, \cite{H2}.
Consider the map $p_1=\mathcal P_1\circ q$, then the main tool is an
extension theorem by Borel, which implies that $p_1$ can be extended
to a holomorphic map of $V_{ss}'$ into $\mathcal M^*$ if
$V_{ss}'\backslash V_s$ is locally contained in a divisor with
normal crossing singularities. The following lemma gives this
property for $p_1$.
\begin{Lem}\label{nor}
Let $C$ be a point in $V_{ss}'\backslash V_s$ i.e. $C$ has an
admissible tacnode. Then $V_{ss}'\backslash V_s$ is smooth at $C$ of
codimension $3$.
\end{Lem}
\proof We first assume that $C$ has only one tacnode in $p$. Then,
after a projective transformation, the curve $C$ can be defined by
an equation of degree four of the form (see Appendix):
$$f=\sum_{i+2j\geq 4}a_{ij}x^iy^j,\ a_{02}=1,$$
in affine coordinates $(x,y)$ centered at $p$. The quartic curves in
a neighbourhood $U$ of $C$ in $V$ can be defined by the equations:
\begin{equation}\label{effeuno}f(t)=\sum_{i+2j\geq 4}a_{ij}x^iy^j+\sum{}t_{ij}x^iy^j=0\end{equation}
where $t=\{t_{ij}\}$ is a system of parameters for $U$. We call
$C(t)$ the quartic curve defined by $f(t)$. If $U$ is sufficiently
small, then a curve belongs to $U\cap (V_{ss}'\backslash V_s)$ if
and only if it has a tacnode. We can assume that $C(t)$ has the
tacnode $p(t)$ in the intersection point of the lines:
$$x+s=0,\ y+ux+v=0$$
where $s,u,v$ depend on $t$ and $y+ux+v=0$ is the tangent line at
$p(t)$. Then $f(t)$ can be written in the form:
\begin{equation}\label{effedue}\sum_{i+2j\geq 4}(a_{ij}+b_{ij})(x+s)^i(y+ux+v)^j\end{equation}
for suitable coefficients $b_{ij}$.
We now compare the equations \ref{effeuno} and \ref{effedue} to get relations of the following types:\\
a) for indices $(i,j)$ such that $i+2j\geq 4$:
$$t_{ij}=b_{ij}+\mbox{terms divisible by } s,u \mbox{ or } v;$$
b) for indices $(i,j)$ such that $i+2j<4$:
$$t_{ij}=\mbox{polynomial in } s, u, v, b_{ij}.$$
Notice that equations of type a) can be solved in $b_{ij}$ as
functions of $s, u, v$ and $t_{ij}$, $i+2j\geq 4$.
Three equations of type b), by forgetting higher terms in $u,s,v$ and $b_{ij}$, $i+2j\geq 4$ are given by:\\
b1)  $\ t_{01}=2v+...$,\\
b2)  $\ t_{30}=a_{31}v+4a_{40}s+a_{21}u+...$,\\
b3)  $\ t_{11}=2a_{12}v+2a_{21}s+2u+...$,\\
We prove that these equations are independent. Otherwise we would
get:
$$a_{21}^2-4a_{40}=0.$$
Let $a_{21}=2\alpha$, then $a_{40}=\alpha^2$. This would imply that
$C$ has an equation of the form:
$$f=(y+\alpha x^2)^2+\sum_{i+2j\geq 5}a_{ij}x^iy^j,$$
i.e. $p$ would be an inadmissible tacnode. Thus we can solve b1),
b2), b3) in $u, s, v$ as functions of $t_{01},   t_{30}, t_{11}$ and
$t_{ij}$, $i+2j\geq 4$. Notice that the cardinality of the set of
indices $(i,j)$ with $i+2j<4$ is equal to $6$. Hence we get $3$
independent equations and it can be easily seen that they define a
variety which is smooth at $C$. In case of a quartic curve $C$ with
two tacnodes, similar computations give that $V_{ss}'\backslash V_s$
has two smooth components of codimension $3$ in $C$. By Lemma 6 in
\cite{H2}, we can find two transversal hyperplanes $L_i$, $i=1,2$
such that $V_{ss}'\backslash V_s\subset L_1\cup L_2$ in a
neighbourhood of $C$.
\qed\
\begin{Prop}\label{ss}
The period map $p_1=\mathcal P_1\circ q$ extends to a holomorphic
map:
$$\mathcal P_2:V_{ss}'\lra \mathcal M^*.$$
\end{Prop}
\proof It follows from Lemma \ref{nor} and Borel's extension theorem
(\cite{Bo}, Theorem A and Remark 3.8).
\qed\\

In fact, the image in $\mathcal V$ of the locus of inadmissible plane
quartics is only one point.
\begin{Lem}\label{in}
The quotient morphism $q:V_{ss} \lra \mathcal V$ maps
$V_{ss}\backslash V_{ss}'$ to the point $v_0$ representing the orbit
of double conics.
\end{Lem}
\proof The proof is similar to that of Lemma 11, \cite{H1}. Let $C$
be an inadmissible plane quartic. We can assume $C$ to be defined by
an equation of the form:
$$f=(yz+\alpha x^2)^2+\sum_{i+2j\geq 5}a_{ij}x^iy^jz^{4-i-j}.$$
Let $F:\mathcal C\ra\Delta$ be a one parameter family of semistable
plane quartics with $F^{-1}(0)=C$. The family $F$ can be given by:
$$f(t)=(yz+\alpha x^2)^2+\sum_{i+2j\geq 5}a_{ij}x^iy^jz^{4-i-j}+t\Psi(x,y,z,t),$$
where $t\in\Delta$ and $\Psi$ is an homogeneous form of degree 4 in
$(x,y,z)$ and holomorphic in $t$. Substituting $t^5$ for $t$ and
$(x,ty,t^{-1}z)$ for $(x,y,z)$ we get the new equivalent family:
$$g(t)=(y+\alpha x^2)^2+\sum_{i+2j\geq 5}a_{ij}t^{2i+j-4}x^iy^jz^{4-i-j}+t^5\Psi(x,ty,t^{-1}z,t^5).$$
Notice that $g$ is holomorphic and that $g(0)=(y+\alpha x^2)^2$.
\qed
\begin{Thm}\label{ss2}
The period map $\mathcal P_0$ can be extended to a holomorphic map:
$$\mathcal P_2:\mathcal V\backslash\{v_0\}\lra \mathcal M^*.$$
The subvariety $\mathcal V\backslash (\mathcal V_s\cup \{v_0\})$ is a smooth rational curve parametrizing type II degenerations of $K3$ surfaces and
it is mapped to the boundary of $\mathcal M^*$.
\end{Thm}
\proof The first assertion follows from Proposition \ref{ss} and
Lemma \ref{in}. Let $C$ be a quartic curve in $V_{ss}'\backslash
V_s$ and $\Delta$ be the open unit disc in $\C$. Consider a one
dimensional family of plane quartic curves $F:\mathcal C\ra \Delta$
with smooth general fiber $C_t$, $t\not=0$ and central fiber
$C_0=C$. By \cite{M}, after passing to a ramified covering of
$\Delta$, we can assume that $C_0$ is a quartic curve in a minimal
orbit. Hence, by Lemma \ref{in} we can assume $C_0$ to be the union
of two conics tangent in two points. Notice that the set of plane
quartics which are the union of two tangent conics maps onto a
smooth rational curve in $\mathcal V$ (with $v_0$ in its closure).
Taking the 4:1 cyclic coverings of $\Ps 2$ branched along the curves
$C_t$ we get a family $G:\mathcal X \ra \Delta$ of quartic surfaces
$X_t$ in $\Ps 3$. The central fiber $X_0$ has two singular points of
type $\tilde E_7$ and it is birational to a ruled variety with base
curve of genus 1. In fact, by Theorem 2.4 in \cite{S1}, this
degeneration is a surface of Type II i.e. the monodromy
transformation $N$ satisfies $N^2=0$, $N\not=0$. In
particular, the monodromy transformation has infinite order. This
implies that the class of $C$ in $\mathcal V$ is mapped to the boundary of
$\mathcal M^*$. \qed
\begin{Rem}\label{shah}
In \cite{S1},  Theorem 2.4, J. Shah gives a classification of GIT
semistable quartic surfaces. It can be easily proved that the 4:1
cyclic cover of $\Ps 2$ branched along a plane quartic is a stable
quartic surface if and only if the plane quartic is stable. In this
case, the quartic surface has at most rational double points. Hence,
we only have surfaces of Type I (case A, Type I in Shah's theorem).
In the case of strictly semistable and reduced plane quartics, the
corresponding quartic surface is strictly semistable with isolated
singularities. In particular, we get quartic surfaces of Type II
(case B, Type II, (i) in Shah's theorem). In the case of a double
conic, the 4:1 cyclic cover of the plane branched along this curve
is the union of two quadrics tangent to each other along the
ramification curve. In particular, this surface $X$ has significant
limit singularities (case B, Surfaces with significant limit
singularities, (ii) in Shah's theorem), this means that the order of
the monodromy transformation depends on the family of surfaces
specializing to $X$. Hence the period map can not be extended to the
point $v_0$.
\end{Rem}

\subsection{Applications}\ \\

\noindent\emph{A) Vinberg's surface.} By a result of T. Shioda and H. Inose there is a one-to-one
correspondence between isomorphism classes of singular $K3$ surfaces
and equivalence classes of positive definite even lattices with
respect to the action of $SL_2(\mathbb Z)$ (see Theorem 4,
\cite{SI}). In \cite{V}, \`E.B. Vinberg describes the geometry of the $K3$ surface with transcendental lattice
$$T\cong \langle 2\rangle^{\oplus 2}.$$
In particular, he finds two interesting projective models for this
surface as a 4:1 cover of the plane branched along a stable quartic.
We give a partial version of his result and provide an alternative
proof.
\begin{Prop}[Theorem 2.5, \cite{V}]\label{vin}
The K3 surface with transcendental lattice isomorphic to $T$ has two
non isomorphic $(L_+,\rho)$-polarizations since it can be obtained
both as the 4:1 cyclic cover of $\Ps 2$ branched along a stable
quartic of singular type $(6,0)$ or of type $(0,3)$. The
isomorphism class of its Picard lattice is given by:
$$U\oplus A_1^{\oplus 2}\oplus E_8^{\oplus 2}.$$
\end{Prop}
\proof Let $C$ be a plane quartic with $6$ nodes i.e. the union of
four lines in $\Ps2$. The corresponding $K3$ surface $X_C$ is the
minimal resolution of a quartic surface in $\Ps3$ with $6$ rational
double points of type $A_3$. The surface $X_C$ carries a natural
involution $\tau_{6,0}$ with fixed locus equal to the disjoint union
of $10$ smooth rational curves. By Proposition \ref{stab}, $X_C$ has
maximum Picard number and the Picard lattice of $X_C$ equals the
invariant lattice $L^{6,0}_+$. By \cite{N2}, Theorem 4.2.2 this
lattice has the invariants:
$$\rho_{6,0}=20,\ \ell_{6,0}=2.$$
Hence, by \cite{N1}, Theorem 4.3.2 the isomorphism classes of the
Picard lattice and the transcendental lattice of $X_C$ are given by:
$$Pic(X_C)\cong U\oplus A_1^{\oplus 2}\oplus E_8^{\oplus 2},$$
$$T(X_C)\cong T.$$
By Corollary \ref{irr} (see type $(0,3)$) and \cite{SI}, Theorem 4,
 the $K3$ surface $X_{C'}$ with $C'$ a plane quartic with three ordinary cusps is isomorphic to $X_C$.\qed\\

\noindent\emph{B) The structure of $\Z[i]$-module on the
transcendental lattice.} We now give a geometric interpretation to
the action of the order four isometry $\rho$ on the lattice $L_-$.
\begin{Lem}
The isomorphism class of $L_-$ is given by:
$$L_-\cong \langle 2\rangle^{\oplus 2}\oplus D_4^{\oplus 3}.$$
\end{Lem}
\proof Since $L_-=L^{\perp}_+$ it follows that $r(L_-)=14$ and
$q_{L_+}=-q_{L_-}$ on $A_{L_+}\cong A_{L_-}$. Then Corollary
\ref{irr} gives that $\ell(L_-)=8,\ \delta(L_-)=1$. Hence the proof
follows from \cite{N1}, Theorem 4.3.2 and \cite{N2}, Proposition 3.2.2.\qed\\
\begin{Lem}\label{act}
The lattices $A^{\oplus 2}_1$ and $D_4$ have a unique structure of $\Z[i]$-modules up to isometries. 
For suitable bases in $A^{\oplus 2}_1$ and $D_4$ this action is
given by the matrices $J_2$ and $J^{\oplus 2}_2$ respectively, where
$$J_2=\left(\begin{array}{cc}
0 &1\\
-1 & 0
\end{array}\right).$$
\end{Lem}
\proof The result for $A^{\oplus 2}_1$ follows from easy
computations. We now consider the case of $D_4$. Recall that the
root lattice $D_4$ can be defined as:
$$D_4=\{x\in\Z^4:\ \sum{x_i}\equiv 0\ mod\,2\},$$
with the standard euclidean inner product.
A basis of $D_4$ is given by the roots:
$$e_1-e_2,\ e_2-e_3,\ e_3-e_4,\ e_3+e_4.$$
Let $f:\mathcal C\ra S$ be a family of stable plane quartics over
the interval $I=[0,1]$ with fiber $C(t)$, $t\not=0$ of type
$(1,0)$ and special fiber $C(0)$ of type $(0,1)$. Then, by Ehresmann's fibration theorem, we have a natural isomorphism:
$$H^2(X_{C(t)},\mathbb Z)\cong H^2(X_{C(0)},\mathbb Z),$$
for every $t\in I$. We fix an isomorphism $H^2(X_{C(t)},\mathbb
Z)\cong L_{K3}$. By Proposition \ref{stab}
$$Pic(X_{C(t)}) \cap L_- \cong A^{\oplus 2}_1,$$
$$Pic(X_{C(0)}) \cap L_- \cong D_4.$$
Hence the family $f$ gives an embedding $A^{\oplus 2}_1\subset D_4.$
By \cite{N1}, Proposition 1.15.1  such embedding is unique up to
isometries. Hence we can assume that $r=e_1-e_2$ and
$\rho(r)=e_3-e_4,$ so $<r,\rho(r)>^{\perp}=<e_1+e_2,e_3+e_4>\cong
A^{\oplus 2}_1.$ Notice that $\rho$ preserves both copies of
$A^{\oplus 2}_1$ and, by the result for $A^{\oplus 2}_1$, it acts on
each copy by the matrix $J_2$ with respect to the natural bases.
This gives to $D_4$ the structure of a free $\Z[i]$-module of rank two and a $\Z[i]$-basis is given by:
$$f_1=e_1-e_2,\ f_2=e_1-e_3,$$
since $\rho(f_1)=e_3-e_4$, $\rho(f_2)=e_3+e_1$.
The matrix representing $\rho$ with respect to the basis $f_1,\rho(f_1),f_2,\rho(f_2)$ is $J_2\oplus J_2$.\qed
\begin{Prop}
The isometry $\rho$ on $L_-$ preserves $\langle 2\rangle^{\oplus 2}$ and each copy of $D_4$, its action is thus described by Lemma \ref{act}.
\end{Prop}
\proof
Consider a stable plane quartic $C$ of type $(0,3)$. Each cusp
gives a lattice of type $D_4\subset Pic(X_C)\cap L_-$ (Proposition \ref{stab}). Moreover $T(X_C)\cong \langle 2\rangle^{\oplus
2}$ (Proposition \ref{vin}). Hence we have that:
$$L_-\cong T(X_C)\oplus (Pic(X_C)\cap L_-).$$
By deformation to plane quartics of singular type $(0,2)$, it is
clear that the automorphism $\rho$ preserves each copy of $D_4$ and
the lattice $T(X_C)$.\qed

\section{Hyperelliptic genus three curves}
In the previous section we proved that the period map can be
extended holomorphically to the complement of the point $v_0$
representing double conics in $\mathcal V$. In Remark \ref{shah} we
observed that the period map can not be extended to this point,
since the 4:1 cyclic cover of the plane branched along a double
conic has significant limit singularities. This is connected to the
existence of \emph{hyperelliptic} $(L_+,\rho)$-polarized $K3$
surfaces i.e.\! such that the curves in the linear system defined by
the degree four polarization introduced in section 1.1. are
hyperelliptic (see \cite{D}). In fact, there is a correspondence
between these polarized hyperelliptic $K3$ surfaces and
hyperelliptic genus three curves which is described in detail by S.
Kond\=o in \cite{K} and \cite{Kh}. In this section we recover
Kond\=o's results providing an alternative description for
hyperelliptic $(L_+,\rho)$-polarized $K3$ surfaces.


\subsection{Smooth hyperelliptic curves}\label{hg}
Let $C\subset \Ps 5$ be a smooth hyperelliptic genus three curve
embedded by the bicanonical map and $i$ be the hyperelliptic
involution on $C$. Consider the surface $\Sigma$ in $\Ps 5$ defined
as follows
$$\Sigma=\overline {\bigcup_{p\in C}\langle p,i(p)\rangle},$$
where $\langle\,,\,\rangle$ denotes the line spanned by two points.
\begin{Prop}\label{hyp}
The surface $\Sigma$ is a cone over a rational normal quartic and
$C$ is a quadratic section of $\Sigma$ not passing through the
vertex. Let $\widetilde \Sigma$ be surface obtained by blowing up of
the vertex of $\Sigma$. This is a $4$-th Hirzebruch and the class of
the curve $C$ in $Pic(\widetilde \Sigma)$ is given by:
$$C=2S_{\infty}+8F,$$
where $S_{\infty}$ is the class of the section with $S_{\infty}^2=-4$ and $F$ is the class of a fiber.
\end{Prop}
\proof
A hyperelliptic genus three curve $C$ can be given by an equation of the form
$$y^2=\prod_{i=1}^{8}(x-\lambda_i),$$
for some complex numbers $\lambda_1,\dots, \lambda_8$.
The hyperelliptic involution $i$ can be written as
$$i:C\lra C,\ \ (x,y)\longmapsto(x,-y).$$
and the bicanonical map is given by
$$\phi_{\mid 2K_C\mid}:C\lra \Ps 5,\ \ (x,y)\longmapsto(1, x, x^2, x^3, x^4, y).$$
Hence, if $(z_0,\dots,z_5)$ are coordinates for $\Ps 5$, the surface
$\Sigma$ is a cone with vertex $(0,\dots,0,1)$ over the rational
normal quartic obtained by projecting $C$ on the hyperplane $z_5=0$.
Moreover, since the curve $C$ is a quadratic section of $\Sigma$,
its inverse image in $\widetilde \Sigma$ has the intersection
numbers: $(C,S_{\infty})=0,$ $(C,F)=2.$ This determines the class of
$C$ in $Pic(\widetilde \Sigma)$.
\qed\\

Let $X_C$ be the 4:1 cyclic cover of the rational ruled surface $\widetilde \Sigma$ branched along the reducible curve:
$$C\cup 2S_{\infty}\in\mid 4S_{\infty}+8F\mid.$$
It can be easily proved that $X_C$ is a $K3$ surface.
\begin{Rem}
In \cite{K} and \cite{Kh} S. Kond\=o associates a $K3$ surface to a smooth hyperelliptic genus three curve in the following way.
The hyperelliptic curve $C$ is embedded in $\Ps1\times \Ps1$ as a divisor of bidegree $(4,2)$, in coordinates $((x_0\!:\!x_1\!),(y_0\!:\!y_1\!))$ for $\Ps1\times \Ps 1$:
$$y_0^2\prod_{i=1}^{4}(x_0-\lambda_ix_1)+y_1^2\prod_{i=5}^{8}(x_0-\lambda_ix_1)=0.$$
Let $L_i=\{y_i=0\}$, $i=1,2$ and $S'_C$ be the double cover of $\Ps1\times \Ps1$ branched along the divisor $B=C+L_1+L_2$ of bidegree $(4,4)$. The minimal resolution $X'_C$ of $S'_C$ is a $K3$ surface.

In fact, it easy to realize that this constructions leads to the
same $K3$ surface obtained before, i.e. $X_C=X'_C$. We give a sketch
of the proof.
Let $C$ be a hyperelliptic genus three curve embedded in a cone $\Sigma$ as in Proposition \ref{hyp}. Let $p_1,\dots,p_8$ be the intersection points of $C$ with the rational normal quartic. It is easy to see that the eight fibers $F_i$ of $\widetilde \Sigma$ through the points $p_i$ are tangent to the curve $C$.
Let $\phi:S_C\lra \widetilde \Sigma$ be the double cover branched along $C$.
The inverse image of the section $S_{\infty}$ is the union of two disjoint smooth rational curves $S_1, S_2$, while the inverse images of the fibers $F_i$ split in the union of two smooth rational curves intersecting transversally in a point on the proper transform of $C$:
$$\phi_1^{-1}(F_i)=F_{i1}\cup F_{i2},\ \ \ i=1,\dots,8.$$
This gives $16$ $(-1)$-curves on $S_C$. We can assume that:
$$F_{ij}\cap S_k\not=\emptyset \Leftrightarrow j=k,\ \ \ \ \
i=1,\dots,8;\ j,k=1,2.$$
We still denote by $F$ the inverse image in $S_C$ of a fiber of $\widetilde \Sigma$.
Let $b:S_C\lra R$ be the blowing down of the eight $(-1)$-curves $F_{11},\dots,F_{41},F_{52},\dots,F_{82}.$
The surface $R$ is a rational ruled surface with sections $L_1=b(S_1)$
and $L_2=b(F)$
Notice that ${L_1}^2=(S_1+F_{11}+\dots+F_{41})^2=0$, ${L_2}^2=0.$
Hence it follows that $R\cong \Ps1\times \Ps1$, $S_C=S'_C$ and the $K3$ surface $X_C=X'_C$ is the double cover of $S_C$ branched along the divisor $C+S_1+S_2$.
\end{Rem}

The previous construction associates a polarized $K3$ surface $X_C$ to any smooth hyperelliptic genus three curve $C$.
In \cite{K} and \cite{Kh} Kond\=o proved that this defines a
holomorphic period map
$$\mathcal P_0^h:\mathcal M^h_3\lra \mathcal D_h\subset \mathcal M.$$
Moreover, if we denote with $\mathcal D'$ the discriminant locus in $\mathcal D^h$ (i.e. the union of hyperplane sections of $\mathcal D_h$ defined by roots), then:
\begin{Thm}[Theorem 5.3, \cite{K}, Theorem 3.8, \cite{Kh}]\label{h1}
The period map induces an isomorphism:
$$\mathcal P_0^h:\mathcal M_3^{h}\lra \mathcal D_h\backslash \mathcal D'.$$
\end{Thm}
\subsection{Stable curves}
As stated in the Appendix, there is an isomorphism
$$\mathcal M^h_3\cong \mathcal V^h_0=V^h_0/\!/PGL_2,$$
where $V^h_0$ denotes the space of sets of eight distinct points in
the projective line. A natural compactification for the moduli space
of smooth hyperelliptic genus three curves is given by the
categorical quotient
$$\mathcal V^h=V^h_{ss} /\!/PGL_2.$$
Recall that a collection of eight points is \emph{stable} (\emph{semistable}) if and only if at most three (four) points coincide. This clearly induces a notion of (semi)stability for hyperelliptic genus three curves.
We prove that there is an isomorphism between $\mathcal V^h$ and a GIT quotient of the space of quadratic sections of a cone which preserves the sets of (semi)stable points.
Let $\Sigma\subset \Ps 5$ be a cone over a rational normal quartic.
\begin{Prop}\label{sec}
The space $\mathcal V^h$ is isomorphic to a GIT quotient of the space of quadratic sections of $\Sigma$ not passing through the vertex with respect to the action of the automorphism group of $\Sigma$.
\end{Prop}
\proof Let $\mathcal G$ be the group of automorphisms
of $\Sigma$, by \cite{De1} this can be described as:
$$\mathcal G =G_0\cdot GL_2/\mu_4,$$
where $G_0\cong H^0(\Ps 1,\mathcal O_{\Ps 1}(4))$, $\mu_4=\{\alpha
I:\alpha^4=1\}$ and the action of $GL_2$ on $G_0$ is the natural
one. The action of $\mathcal G$ can be extended to an action on $\Ps
5$ and $\mathcal O_{\Ps 5}(1)$ admits the following $\mathcal
G$-linearization. Let $z$ be the vertex of $\Sigma$ and
$$A=H^0(\Ps 5, \mathcal O_{\Ps 5}(1))\cong H^0(\Sigma,\mathcal O_{\Sigma}(1)),$$
$$A_1=\mbox{ sections of }A\mbox{ vanishing on }z.$$
Notice that, by restricting the proper transforms of sections to
$S_{\infty}$ in the blow up $\widetilde \Sigma$ we have:
$$ A_1\cong H^0(\Ps1,\mathcal O_{\Ps 1}(4)).$$
Let $\{q_0,q_1,\dots, q_5\}$ be a basis of $A$ where
$\{q_1,\dots,q_5\}$ is a basis for $A_1$.
The action of $G_0$ on $A$ is given by matrices of the form:
$$\phi(\b g)=\left( \begin{array}{cc}
1 & \b g \\
0 & I_5
\end{array}\right),$$
where $\b g$ is a 5-length row vector and $I_5$ the $5\times 5$
identity matrix. Notice that $G_0$ fixes the vertex $z$ and
preserves lines through $z$. The action of $GL_2$ is trivial on
$\mathbb Cq_0$ and acts on $A_1$ via the isomorphism $A_1\cong
H^0(\Ps1,\mathcal O_{\Ps 1}(4)).$ In fact, $GL_2$ fixes the section
$q_0=0$ and acts on it as the group of automorphisms of the rational
normal quartic curve. In order to define the GIT quotient we now
consider the spaces:
$$B=H^0(\Ps 5, \mathcal O_{\Ps 5}(2)),$$
$$B_1 =\mbox{ sections of }B\mbox{ vanishing on }z,$$
$$B_2=\mbox{ sections of }B\mbox{ vanishing with order two on }z.$$
As in the case of $A_1$, by restricting to $S_{\infty}$ we have an
identification  :
$$B_2\cong H^0(\Ps1,\mathcal O_{\Ps1}(8)).$$
We have the following $GL_2/\mu_4$ invariant decomposition:
$$B\cong \mathbb Cq_0^2\oplus q_0A_1\oplus B_2.$$
Consider the open subset of $\mid B\mid$:
$$B^*=\mid B\mid- \mid B_1\mid.$$
A point in $B^* $ is represented uniquely by an element of the
form:
$$q_0^2+q_0a+b,$$
where $a\in A_1$ and $b\in B_2$. If $\mathcal A_1$ and $\mathcal
B_2$ are the affine spaces associated to $A_1$ and $B_2$
respectively, then $B^* $ is isomorphic to $\mathcal A_1\times
\mathcal B_2$. Now we show that, up to the action of $G_0$, every
element in $B^*$ can be represented in the form $q_0^2 +b,$ where
$b\in B_2$. Let $\phi\in G_0$ taking $q_0$ to $q_0+a_0$ with $a_0\in
A_1$. Then $\phi$ takes $f=q_0^2+q_0a+b$ to the form:
$$\phi(f)=q_0^2+(a+2a_0)+(a_0^2+a_0a+b).$$
The choice $a_0=-a/2$ gives the statement and defines a map:
$$\mathcal A_1\times \mathcal B_2\lra \mathcal B_2\ \ \
(a,b)\longmapsto b-a^2/4.$$ It can be be proved as in \cite{S1},
Lemma 4.1, that this map defines a geometric quotient by the action
of $G_0$.

Let $G'$ denote the center of $GL_2/\mu_4$ (isomorphic to the
1-dimensional multiplicative group) and $G''\cong PGL_2$ be the image of
$SL_2$ in $GL_2/\mu_4$. It can be easily seen that:
$$GL_2/\mu_4\cong G'\times G''.$$
By taking the quotient for the action of $G'$ we get
$$\mathcal B_2\backslash\{0\}/\!/G'\cong \mid \mathcal O_{\Ps1}(8)\mid.$$
Moreover, the induced action of $G''$ equals the natural action of
$PGL_2$ on $\mid\mathcal O_{\Ps1}(8)\mid$. Hence, by taking the
categorical quotient of the semistable locus  $\mid \mathcal
O_{\Ps1}(8)\mid_{ss}$ for the action of $PGL_2$ we get a moduli
space isomorphic to $\mathcal V^h$.\qed\\



We now give a characterization of stability as follows:
\begin{Prop}\label{cone}
A quadratic section of $\Sigma$ not passing through the vertex
is stable if and only if it has at most ordinary nodes and cusps.
%
\end{Prop}
\proof Let $p\in C$ and $F$ be the class of a fiber in $\widetilde
\Sigma$. We choose a basis $\{u,v\}$ of $H^0(\widetilde{\Sigma},F)$
such that $u$ vanishes at $p$. Let $l_0=\{u=0\}$,
$l_{\infty}=\{v=0\}$, $q_{\infty}\in
H^0(\widetilde{\Sigma},S_{\infty})$ and $q_0\in
H^0(\widetilde{\Sigma}, 4F+S_{\infty})$. We consider the affine
coordinates $x=u/v$ and $y=q_0/q_{\infty}$. The divisor $C$ is
defined by an equation of the form $f=q_0^2+b=0$, where $b\in
H^0(\widetilde{\Sigma},F)$ and has multiplicity $\leq 4$ at $p$. In
the affine set $\widetilde{\Sigma}-S_{\infty}-l_{\infty}$ the
divisor $C$ is defined by an equation of the form $y^2+p_8(x)=0$
where $p_8$ is a polynomial of degree $8$ which is not divisible by
$x^5$. Notice that $C$ is reduced and has at most double points. If
$C$ is stable then $x^4$ doesn't divide $p_8$, hence $C$ has at most
a node or an ordinary  cusp in $p$.
\qed\\

The minimal resolution of the 4:1 cyclic cover of
$\widetilde{\Sigma}$ branched along the curve $C\cup 2S_{\infty}$,
where $C$ is a stable point in $V^h$ is an $(L_+,\rho)$-polarized K3
surface.
An application of Brieskorn's result as in the proof of Proposition \ref{stab} gives
\begin{Lem}
The period map $\mathcal P_0^h$ extends to a holomorphic map:
$$\mathcal P_1^h:\mathcal V^h_s\lra \mathcal D_h.$$
\end{Lem}
As in the case of stable plane quartics, a natural stratification
for $\mathcal D_h$ is induced by the number of nodes and cusps of a
stable hyperelliptic curve. The Picard lattice of the generic $K3$
surface in each stratum is given in \cite{Kh}, \S4.5.

\subsection{Semistable curves}
Let $C$ be a stricly semistable hyperelliptic genus three curve in a
minimal orbit. With the notation of the proof of Proposition
\ref{cone}, a local equation for the embedding of the curve $C$ in
$\widetilde \Sigma$ is given by:
$$y^2-p_8(x)=0,$$
where $p_8(x)=(x-a)^4(x-b)^4$, $a,b\in \C$. In fact, note that the
locus $\mathcal V^h\backslash \mathcal V^h_s$ is just one point. The
curve $C$ in $\widetilde \Sigma$ has two tacnodes, this implies that
the 4:1 cyclic cover $X_C$ of $\widetilde\Sigma$ branched along
$C\cup 2S_{\infty}$ has two elliptic double points. As in the case
of semistable plane quartics, it can be proved that the period map
can be extended to $\mathcal V^h$ by considering the Baily-Borel
compactification $\mathcal D^*_h$ of the ball quotient $\mathcal
D_h$. In fact, the boundary of $\mathcal D^*_h$ is given by one
point (the \emph{cusp}) and we have the following result.
\begin{Thm}[S. Kond\=o, Theorem 4.7, \cite{Kh}]\label{hiso}
The period map can be extended to an isomorphism:
$$\mathcal P^h:\mathcal V^h\lra \mathcal D^*_h.$$
The point $\mathcal V^h\backslash \mathcal V^h_s$ is mapped to the
cusp.
\end{Thm}
\section{A period map for genus three curves}
In this section we construct a blow-up of $\mathcal V$ in $v_0$ such
that the exceptional divisor is isomorphic to $\mathcal V^h$ and we
prove that the period map $\mathcal P_2$ can be extended to this
variety. In fact, we prove that this extension coincides with the
period map $\mathcal P^h$ on the exceptional divisor.
\subsection{Blow-up}
Let $T$ be a conic in $\Ps 2$ and $2T$ be the corresponding double
conic. We denote by $G_3(2T)$  the orbit of $2T$ by the action of
$PGL_3$.
\begin{Prop}\label{bl}
Let $\widetilde V_{ss}$ be the blow-up of $V_{ss}$ along $G_3(2T)$.
For a proper choice of the $PGL_3$-linearization on $\widetilde V$,
the fibre of
$$b:\widetilde V_{ss}\lra V_{ss}$$
over $2T$ is the semistable locus $V^h_{ss}$
(with respect to the action of $PGL_2$).
The exceptional divisor $\mathcal E$ of the induced blowing-up
$$\widetilde V_{ss}/\!/PGL_3\lra \mathcal V=V_{ss}/\!/PGL_3$$
is isomorphic to $\mathcal V^h$.
\end{Prop}
\proof The result follows from a theorem of F. Kirwan (\S 7,
\cite{Ki}) if we prove that the normal space at $G_3(2T)$ in a point
$2T$ is isomorphic to $H^0(\mathcal O_{\Ps1}(8))$ and that, under
this isomorphism, the isotropy group of $2T$ acts on it as $PGL_2$.

Let $q\in H^0(\Ps2, \mathcal O_{\Ps 2}(2))$ be a section vanishing
on $T$. The isotropy group of $T$ in $PGL_3$ can be easily
identified with $PGL_2$. There exists a unique $GL_2$-invariant
decomposition:
$$H^0(\Ps2, \mathcal O_{\Ps 2}(2))\cong \mathbb C q\oplus \Theta,$$ where $\Theta\cong H^0(\Ps 1,
\mathcal O_{\Ps1} (4))$. Moreover:
$$H^0(\Ps 2, \mathcal O_{\Ps 2}(4))\cong \mathbb Cq^2\oplus q\Theta\oplus \Psi,$$
where $\Psi\cong H^0(\Ps 1,\mathcal O_{\Ps1}(8))$.
The normal space at $G_3(2T)$ in $2T$ is isomorphic to $\Psi$.
Besides, up to the isomorphism of $\Psi$ with $H^0(\Ps1, \mathcal
O_{\Ps 1}(8))$, the action of $PGL_2$ is the canonical one.
Hence, by the result of F. Kirwan, the action of $PGL_3$ can be lifted to $\widetilde V$ in such a way that the exceptional divisor over $v_0$ is isomorphic to the universal categorical quotient of $V^h_{ss}$ by the action of $PGL_2$. \qed\\

We denote with $\widetilde{\mathcal V}$ the blowing-up of $\mathcal V$ in $v_0$ given in Proposition \ref{bl}:
$$\widetilde{\mathcal V}= \widetilde V_{ss}/\!/PGL_3.$$
Let $\widetilde{\mathcal V}_s$ be the subvariety corresponding to the stable
locus in $\widetilde V_{ss}$.

\subsection{Final extension}
We can now associate the isomorphism class of a genus three curve to
each point of $\widetilde{\mathcal V}$. If the point is not on the
exceptional divisor, then the curve can be embedded in $\Ps 2$ as a
semistable plane quartic. Otherwise, it represents a semistable
hyperelliptic genus three curve and can be embedded in the cone
$\Sigma$ as a quadratic section not passing through the vertex.
Moreover, we defined two period maps:
$$\mathcal P_2:\widetilde{\mathcal V}\backslash \mathcal E\lra \mathcal M^*\backslash \mathcal
D_h,$$
$$ \mathcal P^h:\mathcal E\lra \mathcal D^*_h.$$
We now prove that $\mathcal P_2$ and $\mathcal P^h$ give a global
holomorphic period map on $\widetilde{\mathcal V}$.
For similar results and methods see \cite{S} and \cite{S2}.

We start considering the case of a one parameter family of smooth
quartic curves degenerating to a double conic $2T$. Let $X$ be the
4:1 cyclic cover of the Hirzebruch surface $\widetilde \Sigma$
branched along $C+2S_{\infty}$, where $C$ is a quadratic section of
the cone $\Sigma$ not passing through the vertex. We define the
\emph{4:1 cyclic cover of $\Sigma$ branched along $C$ and the
vertex} to be the surface obtained by contracting the inverse images
of $S_{\infty}$ in $X$.
\begin{Prop}\label{pa1}
Let $\mathcal C\lra \Delta$ be a family of plane quartics over the
unit complex disk $\Delta$ intersecting transversally the orbit of
double conics in a smooth
double conic over $0\in\Delta$.\\
Let $F:\mathcal X\lra \Delta$ be the family of quartic surfaces in $\Ps 3$ giving the 4:1 cyclic cover of $\Ps 2$ branched along $\mathcal C$.\\ Then there exist a cover $\xi:\Delta\lra\Delta$  ramified over $0\in \Delta$ and a family $F':\mathcal Y\lra \Delta$ such that:\\
i) the fibers of $F'$ and $\xi^*F$ are isomorphic over $t\not=0$;\\
ii) the central fiber $Y_0$ of $F'$ is the 4:1 cover of a cone
$\Sigma\subset \Ps 5$ over a rational normal quartic, branched
along the vertex and a quadratic section not passing through the
vertex.
\end{Prop}
\proof
Let $(x_0,x_1,x_2)$ be coordinates for $\Ps 2$ and $t\in\Delta$. The
equation of the family $\mathcal C$ can be written in the form:
$$\mathcal C:\ \ q^2(x_i)+t\phi(t,x_i)=0,$$
where $\phi(t)\in H^0(\Ps2, \mathcal O_{\Ps 2}(4))$, $\phi(0)$ is
not divisible by $q$ (from the transversality condition) and $q\in
H^0(\Ps2, \mathcal O_{\Ps 2}(2))$ defines a smooth conic $T$. Let
$\xi$ be the base change of order two $t\mapsto t^2$ on $\Delta$ and
call $\mathcal C'$ the family of curves over $\Delta $ obtained as
pull-back of $\mathcal C$ by $\xi$. The 4:1 cyclic cover $\mathcal
X'$ of $\Delta\times \Ps2$ branched along the family $\mathcal C'$
has a natural embedding in $\Delta\times\Ps 3$, in coordinates
$(t;x_0,x_1,x_2,w)$:
$$\mathcal X':\ \ w^4=q^2(x_i)+t^2\phi(t^2,x_i).$$
Consider the embedding:
$$\Delta\times \Ps 3\longrightarrow \Delta\times\Ps 9$$
given by the identity on the first factor and by the $2$-nd Veronese
embedding on the second one. We consider the coordinates on $\Ps 9$:
$$z_{ij}=x_ix_j,\ s=w^2,\ y_k=tx_k,\ \  \ 0\leq i,j,k\leq 2.$$
The Veronese embedding $V(\Ps 3)$ of $\Ps 3$ is given by the
equations:
$$\left\{ \begin{array}{lll}
z_{ij}z_{kl}-z_{ik}z_{jl}&=&0,\\
sz_{ij}-y_iy_j&=&0,\ \ 0 \leq i,j \leq 2.
\end{array}\right.$$
Then the embedding of $\mathcal X'$ in $\Delta\times\Ps 9$ is the
intersection of $V(\Ps 3)$ with the quadratic section:
$$s^2=q^2(z_{ij})+t^2\phi(z_{ij},t^2),$$
where we think $q\in H^0(\Ps 9, \mathcal O_{\Ps 9}(1))$ and
$\phi(t^2)\in H^0(\Ps 9, \mathcal O_{\Ps 9}(2))$. We still denote by
$T$ the hyperplane section defined by $q$ in $\Ps 9$ and consider
the blowing-up of $\Delta\times \Ps 9$ along $\{0\}\times T$. This
gives a family with general fiber isomorphic to $\Ps 9$ and central
fiber given by the union of a copy of $\Ps 9$ and the exceptional
divisor:
$$E\cong\mathbb P(\mathcal O_T\oplus \mathcal O_T(1)).$$
By Grauert's contraction criterion (see \cite{Gr}), the $\Ps 9$
component in the central fiber can be contracted to a point. The
proper transform $\mathcal Y$ of $\mathcal X'$ in this new variety
is the intersection of the cone over $V(\Ps 3)$ with the quadratic
sections:
$$\left\{\begin{array}{ll}
s^2=\epsilon^2+t^2\phi,\\
q-\epsilon t=0.
\end{array}\right.$$
Let $\mathcal S$ be the projection of $\mathcal Y$ on the linear
subspace defined by:
$$y_i=0,\ \ 0\leq i\leq 2.$$
The surface $S_t$ is defined by:
$$\left\{\begin{array}{ll}
s^2=\epsilon^2+t^2\phi,\\
z_{ij}z_{kl}-z_{ik}z_{jl}=0,\ \ 0\leq i,j\leq 2,\\
q-\epsilon t=0,
\end{array}\right.$$
We denote by $\pi_{2,t}:\mathcal Y\lra\mathcal S$ the projection.\\
Consider now the projection $\mathcal W$ of $\mathcal S$ on the
subspace defined by:
$$s=y_i=0,\ \ 0\leq i \leq 2.$$
The surface $W_t$ is defined by:
$$\left\{\begin{array}{ll}
z_{ij}z_{kl}-z_{ik}z_{jl}=0,\ \ 0\leq i,j\leq 2,\\
q-\epsilon t=0.
\end{array}\right.$$
We denote by $\pi_{1,t}:\mathcal S\lra \mathcal W$ the projection
and  $\pi_t=\pi_{2,t}\circ\pi_{1,t}$.
Then we have the following two cases:\\
i) $t\not=0$:\\
The surface $W_t$ is a linear section of the cone over the Veronese
surface not passing through the vertex, hence it is isomorphic to
the Veronese surface. The projection $\pi_{1,t}$ is the double cover
of $W_t$ branched along the section $B_t$ defined by:
$$\epsilon^2+\phi(z_{ij},t^2)=0.$$
In particular $S_t$ is a Del Pezzo surface of degree two. The
projection $\pi_{2,t}$ is the double cover of $S_t$ branched along
$\pi_{1,t}^{-1}(B_t)$. Hence $\pi_t$ is the 4:1 cyclic cover of
$W_t$ branched along $B_t$. The surface $Y_t$ is isomorphic to the
quartic surface $X_t$.
\\
ii) $t=0$:\\
The central fiber $W_0$ is a cone over the rational normal quartic
given by the intersection of the Veronese surface with the
hyperplane $q=0$. Let $\tilde W_0$ be the $4$-th Hirzebruch surface
obtained by blowing up the vertex of $W_0$, $S_{\infty}$ be the
exceptional divisor and $\tilde B_0$ be the proper transform of
$B_0$. The projection $\pi_{1,0}$ is the double cover of $W_0$
branched along the section $B_0$:
$$\epsilon^2+\phi(z_{ij},0)=0.$$
In particular $S_0$ has two singular points $P_{\pm}$ over the
vertex of $W_0$. The projection $\pi_{2,0}$ is the double cover of
$S_0$ branched along $\pi_{1,0}^{-1}(B_0)$ and the points $P_{\pm}$.
Hence $\pi_0$ is the 4:1 cyclic cover of $W_0$ branched along $B_0$
and the vertex of the cone i.e. the blowing up $\tilde Y_0$ of $Y_0$
in the points $\pi_{2,0}^{-1}(P_{\pm})$ is the 4:1 cyclic cover of
$\tilde W_0$ branched along $\tilde B_0+2S_{\infty}$ .\qed

\begin{Cor}
Let $\mathcal F$ be a one-parameter family in $\widetilde{\mathcal
V}$ intersecting transversally the exceptional divisor $\mathcal E$.
Then the period maps $\mathcal P_2$ and $\mathcal P^h$ glue
holomorphically on $\mathcal F$.
\end{Cor}
\proof By \cite{S2}, Proposition 2.1, up to base change, there
exists a family $\widetilde {\mathcal F}$ in $\widetilde{V}_{ss}$
representing $\mathcal F$ (such that the central fiber belongs to a
minimal orbit). Let $w$ be the intersection of $\widetilde{\mathcal
F}$ with the exceptional divisor. Let $\mathcal C$ be the projection
of $\widetilde{\mathcal F}$ to $V_{ss}$. This is a family of plane
quartics $\mathcal C$ with central fiber equal to a double conic
$2T$. By Proposition \ref{pa1}, after a base change of order two, we
can associate to $\mathcal C$ a family of surfaces such that the
general fiber is isomorphic to the 4:1 cyclic cover of $\Ps 2$
branched along $C_t$ and the central fiber is the 4:1 cyclic cover
of a cone $\Sigma$ branched along a quadratic section $B_0$ and the
vertex. In fact, it follows from the proof of Proposition \ref{pa1}
that the quadratic section $B_0$ is the section corresponding to $w$
as given in Proposition \ref{sec}. Since the period map is invariant
by the action of the Galois group of the double cover, by taking the
minimal resolution of surfaces in the family, we get a period map
$\mathcal P_{\mathcal F}:\mathcal F\ra \mathcal M^*$ which is the
gluing of $\mathcal P_2$ and $\mathcal P^h$.
\qed\\

The existence of a global extension follows from the following
version of Hartogs' Theorem (see \cite{Ho}):
\begin{Lem}[Theorem 2.2.8, \cite{Ho}]\label{hart}
Let $f:U\ra \C$ be a function defined in the open set
$U\subset\C^n$. Assume that $f$ is holomorphic in each variable
$z_j$ when the other variables $z_k$, $k\not=j$ are fixed. Then $f$
is holomorphic in $U$.
\end{Lem}
\begin{Thm}\label{end} The period
map $\mathcal P_2$ can be extended holomorphically to a period map:
$$\mathcal P:\widetilde{\mathcal V}\lra \mathcal M^*$$
such that its restriction to the exceptional divisor $\mathcal E$
coincides with the period map $\mathcal P^h$.
Moreover:\\
i) the locus $\widetilde{\mathcal V}\backslash \widetilde{\mathcal V}_{s}$ is a smooth rational curve mapped to the boundary of $\mathcal M^*$;\\
ii) $\mathcal P$ induces an isomorphism $\mathcal
P_{\mid\widetilde{\mathcal V}_s}:\widetilde{\mathcal
V}_s\lra\mathcal M.$
\end{Thm}
\proof Let $\mathcal P$ be the gluing of $\mathcal P_2$ and
$\mathcal P^h$ on $\widetilde{\mathcal V}$.
The map $\mathcal P$ is holomorphic on $\widetilde{\mathcal
V}\backslash\mathcal E$ (by Theorem \ref{ss}) and on $\mathcal E$
(by Theorem \ref{hiso}). Moreover, by Proposition \ref{pa1}, it is
holomorphic on the generic one-dimensional family intersecting
$\mathcal E$ transversally. Hence, by Lemma \ref{hart}, the period
map $\mathcal P$ is holomorphic. Since $\widetilde{\mathcal V}$ is a
compact variety (see \cite{M}) we also have that $\mathcal P$ is
surjective.

Let $B$ be the rational curve given by the image of strictly
semistable quartics in $\mathcal V-\{v_0\}$ and $b$ be the point in
$\mathcal E$ corresponding to effective divisors of the form
$4p_1+4p_2$ in $\mid \mathcal O_{\Ps 1}(8)\mid_{ss} $. We prove that
the point $b$ lies in the closure of the curve $B$. In fact, a point
in $B$ can be given by:
$$(q+tz^2)(q-tz^2)=q^2-t^2z^4$$
where $q=xy-z^2$ and $t\in\C$. From the proof of Proposition
\ref{pa1}, it follows that the corresponding quadratic section on
the cone is given by the union of two hyperplane sections:
$$q^2-z^4=(q+z^2)(q-z^2)=0.$$
The equation $z^4=0$ intersects the conic $q=0$ in two points with
multiplicity $4$, hence this quadratic section corresponds to the
point $b$.

Moreover, note that the injectivity of the period map on $\mathcal
V_s$ follows as in the generic case. Hence assertions $i)$ and ii)
follow from the surjectivity of the period map and theorems
\ref{ss2} and \ref{hiso}.
\qed\\

Two easy corollaries are the following:
\begin{Cor}
The boundary of the Baily-Borel compactification $\mathcal M^*$ is
given by a unique point (the cusp).
\end{Cor}
\proof This follows from Lemma \ref{bou} and  Theorem \ref{end}.
\qed
\begin{Cor}
Every $(L^+,\rho)$-polarized K3 surface $X$ carries an order four
automorphism $\phi$ such that $\phi^*=\pm iI$ on $H^{2,0}(X)$.
\end{Cor}
\proof By Theorem \ref{end} an $(L_+,\rho)$-polarized $K3$ surface
is either the 4:1 cyclic cover of $\Ps 2$ branched along a stable
quartic or the 4:1 cyclic cover of the cone $\Sigma$ branched along
a quadratic section. In the first case, the $K3$ surface carries an
order four covering automorphism $\sigma_q$ such that the quotient
of $X$ by the involution $\tau_q=\sigma^2_q$ is a rational surface
(the blow up of a Del Pezzo surface of degree two). Hence,
$\tau^*_q$ acts as minus the identity on the transcendental lattice
$T(X)$, as in the generic case. Similarly, in the second case we
have an order four covering automorphism $\sigma_h$ and the quotient
of $X$ by the action of the involution $\tau_h=\sigma_h^2$ is a
rational surface (a blow up of $\Ps1\times\Ps1$).
\qed

\begin{Rem}
In \cite{Lo} E. Looijenga defines a special compactification of the
moduli space $\mathcal M$ which is isomorphic to the GIT
compactification $\mathcal V=V_{ss}/\!/PGL_3$ of the space of plane
quartics. This is obtained by considering a small blowing up and a
blowing down of the Baily-Borel compactification. These correspond
to the blowing up of the cusp and the blowing down of the mirror
$\mathcal D_h$.

Moreover, the same author proved in \cite{Lo2} that the algebra of
invariants of plane quartics can be identified with an algebra of
meromorphic automorphic forms on the complex 6-ball.
\end{Rem}
\appendix
\section{GIT of genus three curves}
The moduli space of genus three curves can be constructed in
different ways, by using geometric invariant theory. For $n\geq 3$ a
genus three curve $C$ can be embedded in $\Ps N$, $N=2(2n-1)-1$ by
the $n$-canonical morphism i.e. the map associated to $\mid
nK_C\mid$. Thus, the moduli space $\mathcal M_3$ can be obtained as
a quotient of either an Hilbert scheme or a Chow variety by the
action of $PGL_{N+1}$. These quotients are compact varieties, hence
they give natural compactification for $\mathcal M_3$. D. Mumford
and D. Gieseker proved that, for $n\geq 5$ the associated
compactification is the moduli space of \emph{Mumford-Deligne stable
curves} (see \cite{M2}). If $n=3$, the associated compactification
is a moduli space for \emph{pseudo-stable curves} described by D.
Schubert in \cite{Sch}. For $n=1$ the GIT construction gives a
compactification for the space of plane quartic curves.
\subsection{Plane quartics}
Let $V=\mid \mathcal O_{\Ps2}(4)\mid\cong\mathbb P^{14}$ be the space of plane quartic curves. We consider the categorical quotient of $V$ with respect to the canonical action of $PGL_3$. We adopt the following notation:
\begin{itemize}
\item[] $V_{ss}=$ set of semistable points in $V$
\item[] $V_{s}=$ set of stable points in $V$.
\item[] $V_0=$ set of smooth plane quartic curves.
\end{itemize}
These loci can be described as follows (see \cite{M}):
\begin{Lem}\label{gt}
Let $C\in V$ be a plane quartic, then:
\item[i)] $C\in V\backslash V_{ss}$ if and only if it is a cubic with an inflectional tangent line or has triple points.
\item[ii)] $C\in V_s$ if and only if it has at most ordinary nodes and cusps.
\item[iii)] $C\in V_{ss}\backslash V_s$  if it is a double conic or it has tacnodes. Moreover, $C$ belongs to the minimal orbit if and only if it is either a double conic or the union of two tangent conics (where at least one is smooth).
\end{Lem}
We consider the equation of a plane quartic in affine coordinates
$(x,y)$ for $\Ps 2$:
$$f=\sum_{i+j\leq 4}a_{ij}x^iy^j.$$
\begin{lemma}
A quartic plane curve $C$  is unstable if and only if, up to a change of coordinates, its equation satisifies one of the following conditions:\\
$I1)\ \ a_{ij}=0\ \ for\ i+j\leq 2,$\\
$I2)\ \ a_{ij}=0\ \ for\ i+2j<4,\ a_{21}=a_{40}=0.$\\
It is strictly semistable if and only if:\\
$SS)\ \ a_{ij}=0\ \ for\ \ \ i+2j<4 ,\ a_{02}\not=0 \mbox{ and }\
a_{21}\not=0 \mbox{ or } a_{40}\not=0.$
\end{lemma}
Notice that condition $I1)$ means that $C$ has a triple point in the
origin, condition $I2)$ corresponds to the case of the union of a
cubic and its inflectional line. The condition $SS)$ means that $C$
has a tacnode.

A double point $p$ on a strictly semistable plane quartic $C$ is
\emph{inadmissible} if there is a projective transformation which
carries $p$ to the point $(0,0,1)$ such that $C$ has an affine
equation of the form:
$$f=(y+\alpha x^2)^2+\sum_{i+2j> 4}a_{ij}x^iy^j.$$
Notice that any point on a double conic is inadmissible and not all
tacnodes are inadmissible.

By \cite{M} there exists a universal categorical quotient
$$q:V_{ss}\lra \mathcal V=V_{ss}/\!/PGL_3,$$
where $\mathcal V$ is a projective variety and $q$ is a
$PGL_3$-invariant surjective morphism
such that its restriction to $V_s$ is a geometric quotient. We
define $\mathcal V_0=q(V_0)$ and $\mathcal V_s=q(V_{s})$. Then there
is a natural isomorphism (\cite{Ve}):
$$V_{0}/\!/PGL_3\cong \mathcal M_3\backslash \mathcal M_3^h,$$
where $\mathcal M_3^h$ denotes the hyperelliptic locus.

The relationship between $\mathcal V$, the Deligne-Mumford compactification and Schubert's compactification are described by D. Mumford in \cite{M2} and by A. Yukie in \cite{Y}.
In particular A. Yukie constructs a blow-up of $\mathcal V$ which maps to the Deligne-Mumford compactification.
\subsection{Hyperelliptic genus three curves}
The subset of smooth hyperelliptic curves is an irreducible divisor $\mathcal M_3^h$ in $\mathcal M_3$.
The canonical morphism gives a natural correspondence between this
locus and the set of collections of eight (unordered) distinct
points in $\Ps 1$ up to projectivities.

Consider the space $V^h=\mid \mathcal O_{\Ps1}(8)\mid $ of effective
divisors of degree eight in $\Ps1$ with the action of the group
$PGL_2$ and let $V^h_0$ be the set corresponding to collections of
eight distinct points. A natural moduli space for smooth
hyperelliptic genus three curves is then given by the geometric
quotient:
$$\mathcal M_3^h\cong V^h_0/\!/PGL_2.$$
We use the notation:
\begin{itemize}
\item[] $V^h_{ss}=$ set of semistable points in $V^h$
\item[] $V^h_{s}=$ set of stable points in $V^h$.
\end{itemize}
These loci can be described as follows (see \cite{M}):
\begin{Lem}
Let $f\in V^h$. Then:
\begin{itemize}
\item[1)] $f\in V^h_s$ iff it has no point of multiplicity
$\geq 4$.
\item[2)] $f\in V^h_{ss}\backslash V^h_s$ if it has points of multiplicity $4$. Moreover, it
belongs to a minimal orbit iff it has two distinct points of
multiplicity $4$.
\end{itemize}
\end{Lem}
A natural compactification for $\mathcal M^h_3$ is the categorical
quotient:
$$\mathcal V^h=V^h_{ss}/\!/PGL_2.$$
We denote with $\mathcal V^h_s$ the locus corresponding to stable
points, note that $\mathcal V^h\backslash \mathcal V^h_s$ is one
point.


\end{document}